\documentclass[ reqno]{amsart}
\usepackage{graphicx} % Required for inserting images

\title[Geometric Aspects of $C^*$-Extreme Points II]{Geometric Aspects of $C^*$-Extreme Points II}
\author{Neha Hotwani}
\address{Department of Mathematics\\
Shiv Nadar Institution of Eminence, Gautam Buddha Nagar-201314, India}
\email{nehahotwani19@gmail.com, neha.hotwani@snu.edu.in}

\usepackage{amsmath,amsthm, amsfonts, amssymb, setspace, color, enumerate, bbold}
      \newtheorem{theorem}{Theorem}[section]

      \newtheorem{remark}[theorem]{Remark}
      \newtheorem{lemma}[theorem]{Lemma}
      
      \newtheorem{proposition}[theorem]{Proposition}

      \def\cK{\mathcal K}

      \def\cT{\mathcal T}
      \def\cA{\mathcal A}

      \def\cQ{\mathcal Q}
      \def\cK{\mathcal K}

      \def\bb1{\mathbb 1}

\usepackage{ulem}
\usepackage{titlesec}
\usepackage{xstring}

\titleformat{\section}
  {\normalfont\Large\bfseries\centering} % format
  {\thesection}{1em}{}                   % numbering + spacing

\renewcommand{\thesection}{\arabic{section}}
\newcommand{\A}{\mathcal A}

\usepackage{ulem}
\date{}

\usepackage[dvipsnames]{xcolor}

\usepackage[pagebackref,colorlinks,linkcolor=Maroon,citecolor=MidnightBlue,urlcolor=NavyBlue,hypertexnames=true]{hyperref}

    \newtheorem{corollary}[theorem]{Corollary}

\begin{document}
\begin{abstract}
    Building on the results of \cite{HR}, we investigate the relationship between \(C^*\)-extremality and linear extremality in the closed unit balls of general unital \(C^*\)-algebras. We show that, unlike in von Neumann algebras, these two notions of extremality need not coincide in general \(C^*\)-algebras. More precisely, we prove that the unilateral shift is a linear extreme point, but not a \(C^*\)-extreme point, of the closed unit ball of the Toeplitz \(C^*\)-algebra.

We further establish that \(C^*\)-extremality implies linear extremality whenever the \(C^*\)-algebra admits a faithful unital representation on a separable Hilbert space. We also prove this implication for elements admitting a factorization $x=v|x|$, where $v$ is a partial isometry belonging to the algebra, and extend the result to suitable quotient \(C^*\)-algebras, including the Calkin algebra. Finally, we characterize the linear extreme points of the closed unit ball of unital \(C^*\)-subalgebras of \(B(H)\) containing the compact operators.
\end{abstract}
\maketitle

\section{Introduction}
\label{sec:Intro}
% Throughout this article, we only consider unital $ C^*$-algebras that are denoted by $\cA$. We use the notation $\cA_1$ to denote the closed unit ball of $\cA$.

% In a $C^*$-algebra, the following theorem is a well-known characterization of the linear extreme points of the closed unit ball.
% \begin{theorem}
%     \label{t:lin-ext}
%     \cite[Theorem~10.2]{T}
%     Let $\cA$ be a  $C^*$-algebra and $\cA_1$ denote the closed unit ball of $\cA$. Then $x$ is a linear extreme point of $\cA_1$ if and only if $x$ is a partial isometry, and it satisfies 
%    \begin{equation}
%        \label{eq:ext-cond}
%       (\mathbf{1}_\cA-x^*x)\cA (\mathbf{1}_\cA-xx^*) =\{0\}.
%    \end{equation}
% \end{theorem}

The geometry of the closed unit ball of a $C^*$-algebra reflects its algebraic structure in a particularly concrete way, as its linear extreme points can be characterized in terms of partial isometries and their defect projections. The notion of $C^*$-convexity provides a complementary perspective by replacing scalar convex coefficients with elements of the algebra. In this setting, extremality is formulated up to unitary equivalence, and the resulting geometry depends on both the admissible coefficients and the unitaries available in the ambient algebra. Understanding the relationship between linear extremality and $C^*$-extremality is therefore a natural problem at the intersection of Banach space geometry and operator algebra theory. The present article continues the study of this relationship, begun in \cite{HR}, beyond the von Neumann algebra setting.

Throughout this article, $\A$ denotes a complex unital $C^*$-algebra with identity $1_{\A}$, and
$\A_1=\{x\in\A:\|x\|\leq 1\}$
denotes its closed unit ball.
In this article, we only consider unital $C^*$-algebras. We write $\sigma_{\A}(x)$ for the spectrum of $x\in\A$, and $|x|=(x^*x)^{1/2}$. When the ambient algebra is clear, we write $\sigma(x)$ in place of $\sigma_{\A}(x)$. Recall that an element $x\in \cA_1$ is a \textbf{linear extreme point} of $\cA_1$ if whenever $x=(y+z)/2$, with $y, z \in \cA_1$, implies $x=y=z$. 
% We denote the set of linear extreme points of $\A_1$ by $\ext(\A_1)$. 
% Thus $x\in\ext(\A_1)$ precisely when $x=(y+z)/2$, with $y,z\in\A_1$, implies $y=z=x$. 
The following is a well-known characterization of linear extreme points of $\cA_1$, due to Kadison.
\begin{theorem}
[{\cite[Theorem~10.2]{T}}]
\label{thm:lin-ext}
Let $\A$ be a unital $C^*$-algebra and let $x\in\A_1$. Then
$x$ is a linear extreme point of $\A_1$ if and only if
$x$ is a partial isometry satisfying
\begin{equation}
\label{eq:ext-cond}
    (1_{\A}-xx^*)\A(1_{\A}-x^*x)=\{0\}.
\end{equation}
\end{theorem}
In particular, every isometry and every coisometry in $\A$ is a linear extreme point of $\A_1$.

Following Loebl and Paulsen \cite{LP}, an element $x\in \A_1$ is called \textbf{$C^*$-convex combination} of $x_1, \dots, x_n \in \cA_1$ whenever $x$ can be represented as
\[
x= \sum_{j=1}^{n}t_j^*x_jt_j,
\quad
 \sum_{j=1}^{n}t_j^*t_j=1_{\A}.
\]
for some $t_1, \dots, t_n \in \A$.
It is called a \textbf{proper $C^*$-convex combination} when every $t_j$ is invertible in $\A$. An element $x\in\A_1$ is a \textbf{$C^*$-extreme point} of $\cA_1$ if every proper $C^*$-convex decomposition of $x$
is trivial up to unitary equivalence; that is, if 
$x=\sum_{j=1}^{n}t_j^*x_jt_j$ 
% x_j\in\A_1,
% \]
is proper, there exist unitaries $u_1, \dots, u_n \in \cA$ such that 
\[
x_j=u_j^*xu_j, \quad 1\leq j\leq n.
\]
% has the property that, for each $j$, there exists $u_j\in\mathcal U(\A)$ such that $x_j=u_j^*xu_j$.
% We write $\cext(\A_1)$ for the set of these points. 

% Although ordinary convex combinations are special cases of $C^*$-convex combinations, the implication from $C^*$-extremality to linear extremality does not follow directly from the definitions: the former yields unitary equivalence of the summands, whereas the latter requires their equality.
% For the algebra $\mathcal B(H)$ of bounded linear operators on a complex Hilbert space $H$, Hopenwasser, Moore and Paulsen \cite{HMP} characterized the $C^*$-extreme points of the closed unit ball as the isometries and coisometries. More recently, 
Recently, the study of $C^*$-extremality has attracted interest in various
operator-algebraic settings; see
\cite{BBK,BC,BK,GMR,RK,RKS,K,KM,HR2} and the references
therein.
In \cite{HR}, the authors established that $C^*$-extremality and linear extremality coincide in the closed unit ball of every von Neumann algebra,
% In particular, they proved that every isometry and every coisometry is a $C^*$-extreme point of the closed unit ball,
thereby extending the result of \cite{HMP} to arbitrary von Neumann algebras. One of the main results of the present article shows that this conclusion does not extend to general $C^*$-algebras. 
More precisely, we prove the following theorem.
\begin{theorem}
    \label{thm:shift-not-cstar-extreme}
Let $S$ be the unilateral shift on the Hardy space $H^2$ on the unit circle, and let
$\mathcal T=C^*(S)$
be the Toeplitz $C^*$-algebra. Then $S$ is not a $C^*$-extreme point
of the closed unit ball $\mathcal T_1$.
\end{theorem}

% that the unilateral shift operator, although 
A proof of Theorem~\ref{thm:shift-not-cstar-extreme} is given in Section~\ref{sec:main}. Since \(S\) is an isometry in $\cT$, it is a linear extreme point of $\cT_1$. Thus, Theorem~\ref{thm:shift-not-cstar-extreme} shows that an isometry in a unital \(C^*\)-algebra need not be a \(C^*\)-extreme point of its closed unit ball. Consequently, linear extremality does not imply \(C^*\)-extremality in general unital \(C^*\)-algebras. Moreover, since \(S\) remains an isometry in the enveloping von Neumann algebra \(\mathcal T^{**}\), it is a \(C^*\)-extreme point of \((\mathcal T^{**})_1\) by \cite[Theorem~1.2]{HR}. Therefore, an element of a unital \(C^*\)-algebra may be \(C^*\)-extreme in the closed unit ball of its enveloping von Neumann algebra without being \(C^*\)-extreme in the original algebra.

% Related results for algebras with a faithful tracial state and for spaces of continuous algebra-valued functions were obtained in \cite{HRVec26}. 
% These results lead to two questions: how far does the implication from $C^*$-extremality to linear extremality extend to general unital $C^*$-algebras, and does the converse hold outside the von Neumann algebra setting?

Further, we show that if a $C^*$-extreme point of $\cA_1$ is a partial isometry, then it is a linear extreme point of $\cA_1$ (Proposition \ref{p:p.i.+c_ext}). We then prove that whenever \(\cA\) admits a faithful unital $*$-representation on a separable Hilbert space, every \(C^*\)-extreme point of \(\cA_1\) is a partial isometry and hence a linear extreme point (Theorem~\ref{thm:rep-on-sep}). In particular, this implication holds for every separable unital \(C^*\)-algebra. However, the hypothesis does not require \(\cA\) itself to be norm separable. 
Moreover, unlike the results in \cite{HR3}, our theorem does not require a tracial assumption and applies to \(C^*\)-algebras beyond the class of von Neumann algebras. Combining this result with
Theorem~\ref{thm:shift-not-cstar-extreme},
we conclude that every $C^*$-extreme point
of the closed unit ball of the Toeplitz
$C^*$-algebra $\cT$ is a linear extreme point,
but the converse does not hold (Corollary \ref{cor:Toeplitz-alg}).

We also establish an alternative sufficient condition that does not require a separable representation. More precisely, if \(x\in\cA_1\) admits a factorization \(x=v|x|\), where \(v\in\cA\) is a partial isometry, then \(C^*\)-extremality of \(x\) implies linear extremality (Theorem~\ref{thm:A-has-p.d.}). In particular, this implication holds for every unital \(C^*\)-algebra in which each element admits a polar decomposition within the algebra. Furthermore, if every element of a unital \(C^*\)-algebra admits such a factorization, then the same implication holds in each of its unital quotients (Theorem~\ref{cor:ideal}). In particular, we obtain that \(C^*\)-extremality implies linear extremality in the Calkin algebra (Corollary~\ref{cor:cal-alg}).

We conclude the introduction by fixing some notation
that will be used throughout the article.
For a Hilbert space $H$, let $B(H)$ and $\cK(H)$
denote the spaces of bounded linear operators and
compact operators on $H$, respectively.
We write $I_H$ for the identity operator on $H$.
For $T\in B(H)$, the commutant of $T$ is defined by
$\{T\}'=\{R\in B(H):RT=TR\}$.
% If $T\in B(H)$ is self-adjoint, we denote its
% projection-valued spectral measure by $E_T$.
Recall that a $*$-representation
$\pi:\cA\to B(H)$ is faithful if it is injective
and unital if $\pi(1_{\cA})=I_H$.
For a compact Hausdorff space $\Omega$,
$C(\Omega)$ denotes the $C^*$-algebra of
continuous complex-valued functions on $\Omega$,
equipped with the supremum norm $\|\cdot\|_\infty$.

For the Toeplitz $C^*$-algebra, let
\[
\mathbb D=\{z\in\mathbb C:|z|<1\},
\qquad
\mathbb T=\partial\mathbb D.
\]
The disk algebra is defined by
\[
A(\mathbb D)
=
\{g\in C(\overline{\mathbb D}):
g|_{\mathbb D}\text{ is holomorphic}\},
\]
and its elements are identified with their
boundary restrictions to $\mathbb T$.

With respect to normalized arc-length measure
on $\mathbb T$, the Hardy space $H^2(\mathbb T)$
is the closed linear span of
$\{z^n:n\geq0\}$ in $L^2(\mathbb T)$.
We write $H^2$ for $H^2(\mathbb T)$ and denote
by $P_+$ the orthogonal projection from
$L^2(\mathbb T)$ onto $H^2$.

For $f\in L^\infty(\mathbb T)$, the Toeplitz
operator $T_f$ on $H^2$ is defined by
\[
T_fh=P_+(fh),\qquad h\in H^2.
\]
For a bounded analytic function
$g\in H^\infty(\mathbb D)$, we denote by
$M_g$ the multiplication operator on $H^2$,
where $g$ is identified with its almost-everywhere
boundary values on $\mathbb T$.
In particular, the unilateral shift is $S=M_z$.
For standard results concerning Hardy spaces,
Toeplitz operators, and the Toeplitz $C^*$-algebra,
we refer to \cite[Chapters~6 and~7]{D}.
% Standard Hardy-space and Toeplitz-operator facts are used in the form presented in \cite[Chapters~6-7]{D}. These include the commutant identity
% \[
%  \Talg\cap\{S\}'=\{M_g:g\in A(\mathbb D)\},
% \]
% which is central to the proof of Theorem~\ref{thm:shift-not-cstar-extreme}.

\section{Main Results}
\label{sec:main}

% In this section, we prove our main results. 
In \cite{HR}, the authors proved that isometry is always a $C^*$-extreme point of the closed unit ball in a von Neumann algebra by proving the Wold decomposition type theorem as in the case of $B(H)$ in \cite{HMP}. Therefore, one can ask whether an isometry is always a $C^*$-extreme point of the closed unit ball in a $C^*$-algebra. Following the same technique as in \cite[Theorem~1.1]{HMP}, it is equivalent to ask whether two isometries are similar in a $C^*$-algebra; does this imply that they are unitarily equivalent within the $C^*$-algebra? In Theorem \ref{thm:shift-not-cstar-extreme}, we provide an example which shows that this is not true in general. More precisely, we show that in the Toeplitz algebra, the unilateral shift operator is similar to an isometry from the algebra, but it is not unitarily equivalent to that isometry. We begin by proving some lemmas that will be used in Theorem \ref{thm:shift-not-cstar-extreme}. 

\begin{lemma}
\label{lem:disk-algebra}
There exists a function $f\in C(\mathbb T)$ such that
$0<f<1$
and
$f\neq |g|^2$
for every $g\in A(\mathbb D)$.
\end{lemma}

\begin{proof}
By \cite[Corollary 2.3]{BP}, there exists a strictly positive
function $f_0\in C(\mathbb T)$ such that
$f_0\neq |h|$
for every $h\in A(\mathbb D)$. Since $\mathbb T$ is compact, $f_0$ is
bounded. Choose $\lambda>0$ such that
$\lambda \|f_0\|_\infty<1$
and define
$f=(\lambda f_0)^2$.
Then $f\in C(\mathbb T)$ and
$0<f<1$.
If, for some $g\in A(\mathbb D)$, we had
$f=|g|^2$,
then
$\lambda f_0=|g|$,
and consequently
$f_0=\left|\frac{g}{\lambda}\right|$,
where $g/\lambda\in A(\mathbb D)$. This contradicts the choice of
$f_0$. Hence
$f\neq |g|^2$
for every $g\in A(\mathbb D)$.
\end{proof}

The following lemma is perhaps well known; we provide a proof for completeness. We use the notation $\mathbf{1}$ to denote the constant function on $\mathbb T$ that takes the value 1 for all $z\in \mathbb T$. 
\begin{lemma}
    \label{lem:commutant}
    Let $S$ be the unilateral shift operator in $\cT$. Then 
    \[
\mathcal T\cap\{S\}'
=
\{M_g:g\in A(\mathbb D)\}.
\]
\end{lemma}

\begin{proof}
Let $X\in \cT\cap\{S\}'$. Since $XS=SX$ and
$S^*S=I_H$, we have
$X=S^{*n}XS^n$ for $n\geq 1$.
By \cite[Proposition~11.2.19]{E},
\[
\cT=\{T_f+K:f\in C(\mathbb T),
K\in\cK(H^2)\}.
\]
Thus, $X=T_f+K$ for some $f\in C(\mathbb T)$
and $K\in\cK(H^2)$. Consequently,
\[
X=S^{*n}(T_f+K)S^n
  =S^{*n}T_fS^n+S^{*n}KS^n.
\]
From \cite[Proposition~11.2.16]{E}, we have
$S^*T_fS=T_f$. Therefore,
$X=T_f+S^{*n}KS^n$.
Now, $S^{*n}\to 0$ strongly, and the compactness
of $K$ implies that
$\|S^{*n}K\|\to 0$.
As $\|S^n\|=1$, we obtain
$\|S^{*n}KS^n\|
\leq \|S^{*n}K\|\to 0$.
Hence $X=T_f$.
Since $T_fS=ST_f$, we have
\[
0=\langle (T_fS-ST_f)z^n,\mathbf{1}\rangle
  =\widehat f(-n-1),
\qquad n\geq 0.
\]
% Here $\mathbf{1}$ denotes the constant function on $\mathbb T$.
Thus, all negative Fourier coefficients of $f$
vanish. Since \(f\in C(\mathbb T)\), Fejér's theorem (\cite[Theorem~3.1]{Kat}) implies that \(f\) is a uniform limit of analytic polynomials on \(\mathbb T\). Consequently, \(f\) admits an extension belonging to \(A(\mathbb D)\).

Conversely, let $f\in A(\mathbb D)$. By uniform
polynomial approximation, there exists a
sequence of analytic polynomials $(p_n)$ such
that $p_n\to f$ uniformly on
$\overline{\mathbb D}$. Therefore,
\[
\|T_f-p_n(S)\|
\leq \|f-p_n\|_\infty\to 0.
\]
Since $p_n(S)\in\cT\cap\{S\}'$ for every $n$,
and $\cT\cap\{S\}'$ is norm closed, we conclude
that $T_f\in\cT\cap\{S\}'$.
\end{proof}

% We refer to \cite[Chapter~6,7]{D} for standard results related to the Toeplitz algebra.
Now we prove Theorem \ref{thm:shift-not-cstar-extreme}.
\begin{proof}
Let $f\in C(\mathbb T)$ be as in Lemma~\ref{lem:disk-algebra}, and let
$A=T_f$.
% where $T_f$ denotes the Toeplitz operator with symbol $f$.
Since $f$ is continuous and
$0<f<1$,
there exist constants $m,M$ such that
$0<m\leq f\leq M<1$.
Consequently,
$mI_H\leq T_f\leq MI_H<I_H$. 
% Here $I_H$ denotes the identity operator in $\cT$.
Thus both $A=T_f$ and $I_H-A$ are positive and invertible elements of
$\mathcal T$.
Since $S^*AS=A$, one has
% $S^*AS=A$
% and therefore
$S^*(I_H-A)S=I_H-A$.
Let
\[
T_1=A^{1/2},
\qquad
T_2=(I_H-A)^{1/2}.
\]
Then $T_1,T_2\in\mathcal T$ are positive and invertible, and
\[
T_1^*T_1+T_2^*T_2=T_1^2+T_2^2=A+(I_H-A)=I_H.
\]
Define
$V=T_1ST_1^{-1}$,
and
$W=T_2ST_2^{-1}$.
Again using $S^*AS=A$, one obtains
\[
\begin{aligned}
V^*V
&=T_1^{-1}S^*T_1^2ST_1^{-1}\\
&=T_1^{-1}S^*AST_1^{-1}\\
&=T_1^{-1}AT_1^{-1}\\
&=I_H.
\end{aligned}
\]
Thus $V$ is an isometry. Similarly,
$W^*W=I_H$,
so $W$ is also an isometry. In particular,
$V,W\in\mathcal T_1$.
Moreover,
$T_1VT_1
=T_1(T_1ST_1^{-1})T_1
=T_1^2S
=AS$,
and similarly,
$T_2WT_2=(I_H-A)S$.
Therefore
$S=T_1^*VT_1+T_2^*WT_2$.
Since $T_1$ and $T_2$ are invertible and
$T_1^*T_1+T_2^*T_2=I$,
this is a proper $C^*$-convex combination of elements of
$\mathcal T_1$.

Suppose, to the contrary, that $S$ is a $C^*$-extreme point
of $\mathcal T_1$. Then $V$ must be unitarily equivalent to $S$. Hence there exists a unitary
$U\in\mathcal T$ such that
$V=USU^*$.
Since
$V=T_1ST_1^{-1}$,
we have
$T_1ST_1^{-1}=USU^*$,
which gives
$(U^*T_1)S=S(U^*T_1)$.
Put
$C=U^*T_1$.
Then
$C\in\mathcal T\cap\{S\}'$.
% We use the standard identity
% \[
% \mathcal T\cap\{S\}'
% =
% \{M_g:g\in A(\mathbb D)\}.
% \]
Now applying Lemma \ref{lem:commutant}, we obtain
$C=M_g$
for some $g\in A(\mathbb D)$.
On the other hand,
$C^*C
=
T_1UU^*T_1
=
T_1^2
=
A
=
T_f$.
Since $C=M_g$ and $g$ is analytic,
$C^*C
=
M_g^*M_g
=
T_{|g|^2}$.
Thus
$T_f=T_{|g|^2}$.
Equality of Toeplitz operators implies equality of their symbols almost
everywhere. Since both $f$ and $|g|^2$ are continuous on
$\mathbb T$, it follows that
$f=|g|^2$
\text{on }$\mathbb T$.
This contradicts Lemma~\ref{lem:disk-algebra}.
Therefore $S$ is not a $C^*$-extreme point of $\mathcal T_1$.
\end{proof}

The following proposition shows that a \(C^*\)-extreme point of the closed unit ball of a \(C^*\)-algebra is a linear extreme point whenever it is a partial isometry.

\begin{proposition}
\label{p:p.i.+c_ext}
Let $\cA$ be a  $C^*$-algebra and let $x\in \cA_1$ be a partial
isometry. If $x$ is a $C^*$-extreme point of $\cA_1$, then $x$ is a
linear extreme point of $\cA_1$.
\end{proposition}

\begin{proof}
Let
$p=x^*x$
and
$q=xx^*$.
Since $x$ is a partial isometry, $p$ and $q$ are projections. By
Theorem \ref{thm:lin-ext}, it remains to show that $x$ satisfies Equation \eqref{eq:ext-cond}, i.e., to show that
\[
(\mathbf 1_{\cA}-q)\cA(\mathbf 1_{\cA}-p)=\{0\}.
\]
Suppose, to the contrary, that
$(\mathbf 1_{\cA}-q)\cA(\mathbf 1_{\cA}-p)\neq\{0\}$.
Choose
\[
 z\in
(\mathbf 1_{\cA}-q)\cA(\mathbf 1_{\cA}-p)
\]
with $0<\|z\|<1$.
Since $x^*z=z^*x=0$, we have
$(x\pm z)^*(x\pm z)=p+z^*z$.
Moreover, since
$z=z(1_{\cA}-p)$, we obtain
$0\leq z^*z
\leq \|z\|^2(1_{\cA}-p)$.
Consequently,
\[
0\leq p+z^*z
\leq p+\|z\|^2(1_{\cA}-p)
\leq 1_{\cA}.
\]
Thus, $x+z$ and $x-z$ belong to $\cA_1$.
% \[
% 0< z^*z=(\mathbf{1}_\cA -p) z^*z (\mathbf{1}_\cA -p) < \mathbf{1}_\cA -p.
% \]
% In particular, $0< p+z^*z<\mathbf{1}_\cA$.
Now, 
\[
x=\frac{1}{2}\left((x+z)+(x-z)\right).
\]
Equivalently,
\[
x
=
\left(\frac{1}{\sqrt2}\mathbf{1}_{\mathcal A}\right)^*
\left(x+z\right)
\left(\frac{1}{\sqrt2}\mathbf{1}_{\mathcal A}\right)
+
\left(\frac{1}{\sqrt2}\mathbf{1}_{\mathcal A}\right)^*
(x-z)
\left(\frac{1}{\sqrt2}\mathbf{1}_{\mathcal A}\right),
\]
which is a proper $C^*$-convex decomposition of $x$. Since $x$ is a $C^*$-extreme point of $\cA_1$, one has $x$ is unitarily equivalent to $x+z$. 
Consequently,
$(x+z)^*(x+z)=p+z^*z$
is a projection.
Since $pz^*z=z^*zp=0$, we obtain
$(p+z^*z)^2=p+(z^*z)^2$.
Therefore,
$(z^*z)^2=z^*z$.
Thus, $z^*z$ is a projection.
However,
$0<\|z^*z\|=\|z\|^2<1$,
which is impossible.
% Therefore, $p=x^*x$ is unitarily equivalent to $p+z^*z$. Hence $(z^*z)^2=z^*z$ which is a contradiction because $0<z^*z<\mathbf{1}_\cA$.
Hence, $x$ satisfies Equation \eqref{eq:ext-cond}. This completes the proof.
\end{proof}

Next, we show that $ C^*$-extremality implies linear extremality in a $ C^*$-algebra that has a faithful representation on a separable Hilbert space. 
Consequently, the result holds for all separable $C^*$-algebras. The proof uses techniques from \cite[Lemma 2.8]{HR3}.
% For instance, Cuntz algebra, Toeplitz algebra, $C^*(G)$, $G$ countable discrete group have a faithful representation on a separable Hilbert space.

\begin{theorem}
\label{thm:rep-on-sep}
Let $\mathcal A$ be a $C^*$-algebra admitting a faithful
unital $*$-representation
$\pi:\mathcal A\to B(H)$
on a separable Hilbert space $H$. If $x$ is a $C^*$-extreme point
of $\mathcal A_1$, then $x$ is a partial isometry. In particular, it is a linear extreme point of $\cA_1$.
\end{theorem}

\begin{proof}
Let
$p=x^*x$.
Since $x\in \mathcal A_1$, we have $0\leq p\leq 1_{\mathcal A}$.
Suppose, to the contrary, that $p$ is not a projection. Then
$\sigma_{\mathcal A}(p)\cap (0,1)\neq\varnothing$.
% Indeed, if $\sigma_{\mathcal A}(p)\subseteq\{0,1\}$, then the
% continuous functional calculus would imply $p^2=p$.
Choose
$\lambda_0\in \sigma_{\mathcal A}(p)\cap(0,1)$
and set
$T=\pi(p)\in B(H)$.
Since $\pi$ is faithful,
$\sigma_{B(H)}(T)=\sigma_{\mathcal A}(p)$,
and hence $\lambda_0\in\sigma(T)$. Let $E_T$ denote the spectral
measure of $T$.
% We now adapt the spectral perturbation argument from
% \cite[Lemma~2.8]{HR}. 
Consider the continuous function
\[
\Theta:(0,1)\times(0,1)\longrightarrow\mathbb R
\]
defined by
\[
\Theta(\beta,\gamma)
=
\frac{1}{\sqrt{\beta}}-1
-
\left|\sqrt{\frac{\gamma}{\beta}}-1\right|.
\]
Since
\[
\Theta(\lambda_0,\lambda_0)
=
\frac{1}{\sqrt{\lambda_0}}-1>0,
\]
there exists an open interval $U$ containing $\lambda_0$,
with $\overline U\subset(0,1)$, such that
\begin{equation}
\left|\sqrt{\frac{\gamma}{\beta}}-1\right|
<
\frac{1}{\sqrt{\beta}}-1,
\qquad \beta,\gamma\in U.
\label{eq:F-estimate}
\end{equation}
% This is precisely the local estimate used in
% \cite[Lemma~2.8]{HR}.

Since $\lambda_0\in\sigma(T)$, we have $E_T(U)\neq0$. Hence we may
choose a compact interval $J\subset U$ such that
$E_T(J)\neq0$.
As $T$ is self-adjoint and $H$ is separable, $T$ has at most
countably many distinct eigenvalues. Since $U$ is uncountable, we
may choose
$\gamma_0\in U$
which is not an eigenvalue of $T$.
% Following again the construction in
% \cite[Lemma~2.8]{HR}, 
Choose
$\phi\in C([0,1],[0,1])$ such that
\[
\phi|_J=1
\qquad\text{and}\qquad
\phi=0\quad\text{on }[0,1]\setminus U,
\]
and define
\[
\delta(\beta)=
\begin{cases}
\displaystyle
\phi(\beta)\left(\sqrt{\frac{\gamma_0}{\beta}}-1\right),
& \beta\in U,\\[2mm]
0,&\beta\notin U.
\end{cases}
\]
Then $\delta\in C([0,1],\mathbb R)$ and, by
\eqref{eq:F-estimate},
\[
|\delta(\beta)|
\leq
\frac{1}{\sqrt{\beta}}-1,
\qquad 0<\beta\leq1.
\]
Consequently,
$|1\pm \delta(\beta)|
\leq \frac{1}{\sqrt{\beta}}$,
and therefore
\begin{equation}
\beta(1\pm \delta(\beta))^2\leq1,
\qquad 0\leq \beta\leq1.
\label{eq:contraction-estimate}
\end{equation}
By continuous functional calculus,
$\delta(p)\in C^*(\mathbf{1}_{\mathcal A},p)\subseteq\mathcal A$, where $C^*(\mathbf{1}_{\mathcal A},p)$ denotes the smallest $C^*$-algebra generated by $\mathbf{1}_\cA$ and $p$.
Define
$y_{\pm}
=
x\bigl(\mathbf{1}_{\mathcal A}\pm \delta(p)\bigr)$.
Since $\delta(p)$ commutes with $p$, we obtain from
\eqref{eq:contraction-estimate} that
\[
\begin{aligned}
y_{\pm}^*y_{\pm}
&=
\bigl(\mathbf{1}_{\mathcal A}\pm \delta(p)\bigr)
p
\bigl(\mathbf{1}_{\mathcal A}\pm \delta(p)\bigr)\\
&=
p\bigl(\mathbf{1}_{\mathcal A}\pm \delta(p)\bigr)^2
\leq \mathbf{1}_{\mathcal A}.
\end{aligned}
\]
Thus $y_{\pm}\in\mathcal A_1$. Moreover,
\[
x=\frac{y_++y_-}{2},
\]
or, equivalently,
\[
x
=
\left(\frac{1}{\sqrt2}\mathbf{1}_{\mathcal A}\right)^*
y_+
\left(\frac{1}{\sqrt2}\mathbf{1}_{\mathcal A}\right)
+
\left(\frac{1}{\sqrt2}\mathbf{1}_{\mathcal A}\right)^*
y_-
\left(\frac{1}{\sqrt2}\mathbf{1}_{\mathcal A}\right).
\]
This is a proper $C^*$-convex decomposition of $x$. Since $x$ is
$C^*$-extreme, $y_+$ is unitarily equivalent to $x$. Hence there
exists a unitary $u\in\mathcal A$ such that
$y_+=u^*xu$.
It follows that
\begin{equation}
y_+^*y_+
=
u^*pu.
\label{eq:unitary-equivalence}
\end{equation}
Define
\[
\psi(\beta)=\beta(1+\delta(\beta))^2,
\qquad 0\leq \beta\leq1.
\]
Then
$y_+^*y_+=\psi(p)$.
Applying $\pi$ to \eqref{eq:unitary-equivalence}, we obtain
$\psi(T)
=
\pi(u)^*T\pi(u)$,
so $\psi(T)$ and $T$ are unitarily equivalent.
On the other hand,  since $\phi=1$ on $J$, for
$\beta\in J$ we have
\[
\delta(\beta)=\sqrt{\frac{\gamma_0}{\beta}}-1,
\]
and hence
\[
\psi(\beta)
=
\beta\left(\sqrt{\frac{\gamma_0}{\beta}}\right)^2
=
\gamma_0.
\]
Therefore,
\[
\psi(T)E_T(J)=\gamma_0E_T(J).
\]
Since $E_T(J)\neq0$, it follows that $\gamma_0$ is an eigenvalue of
$\psi(T)$. However, $\gamma_0$ was chosen not to be an eigenvalue of $T$.
Thus $\psi(T)$ cannot be unitarily equivalent to $T$, a
contradiction.
Therefore $p=x^*x$ is a projection. Hence $x$ is a partial
isometry. Now applying Proposition \ref{p:p.i.+c_ext}, we obtain that $x$ is indeed a linear extreme point of $\cA_1$.
\end{proof}

% \begin{corollary}
%     \label{cor:rep-on-sep}
%     Let $\cA$ be a unital $C^*$-algebra such that it has a faithful representation on a separable Hilbert space $H$. Let $x$ be a $C^*$-extreme point of $\cA_1$; then $x$ is a linear extreme point of $\cA_1$. 
% \end{corollary}

% \begin{proof}
%    It follows immediately from Theorem \ref{thm:rep-on-sep} and Proposition \ref{p:p.i.+c_ext}. 
% \end{proof}

\begin{corollary}
\label{cor:Toeplitz-alg}
Let $\cT$ be the Toeplitz $C^*$-algebra.
Then every $C^*$-extreme point of $\cT_1$
is a linear extreme point of $\cT_1$,
but the converse does not hold.
\end{corollary}

\begin{proof}
Since $\cT$ admits a faithful unital
$*$-representation on the separable Hilbert
space $H^2(\mathbb T)$,
Theorem~\ref{thm:rep-on-sep} implies that
every $C^*$-extreme point of $\cT_1$
is a linear extreme point.
On the other hand, the unilateral shift $S$
is an isometry in $\cT$ and hence a linear
extreme point of $\cT_1$.
However, by
Theorem~\ref{thm:shift-not-cstar-extreme},
$S$ is not a $C^*$-extreme point of $\cT_1$.
This completes the proof.
\end{proof}

In \cite{HR}, the existence of polar decompositions
in von Neumann algebras played a crucial role in
establishing the relationship between $C^*$-extremality
and linear extremality. The following theorem shows
that, for an element $x\in\cA_1$, the weaker condition
$x=v|x|$, where $v\in\cA$ is a partial isometry,
suffices to ensure that $C^*$-extremality implies
linear extremality.

\begin{theorem}
    \label{thm:A-has-p.d.}
   Let $x\in\cA_1$ admit a factorization
$x=v|x|$ for some partial isometry $v\in\cA$.
If $x$ is a $C^*$-extreme point of $\cA_1$,
then $x$ is a linear extreme point of $\cA_1$.
\end{theorem}

\begin{proof}
   By hypothesis, \(x=v|x|\), where \(v\in\cA\) is a partial isometry. Since \(|x|\) is a positive contraction, it can be expressed as the average of two unitaries \(u_1,u_2\in\cA\). Consequently,

\[
x=\frac{vu_1+vu_2}{2}.
\]

Since \(vu_1\) and \(vu_2\) are partial isometries in \(\cA\), they belong to \(\cA_1\). Moreover, the above expression is a proper \(C^*\)-convex decomposition of \(x\). By the \(C^*\)-extremality of \(x\), there exists a unitary \(u\in\cA\) such that
$x=u^*(vu_1)u$.
Hence \(x\) is a partial isometry in \(\cA\). Applying Proposition~\ref{p:p.i.+c_ext}, we conclude that \(x\) is a linear extreme point of \(\cA_1\).
\end{proof}

\begin{remark}
Theorem~\ref{thm:A-has-p.d.} also provides a simpler
proof of the implication $(i)\implies(ii)$ in
\cite[Theorem~1.2]{HR}. More generally, the theorem
applies to every unital $C^*$-algebra in which each
element admits a polar decomposition within the algebra.
\end{remark}
%  \begin{remark}
% Since every element of a von Neumann algebra admits a polar
% decomposition within the algebra, Theorem~\ref{thm:A-has-p.d.}
% provides a simpler proof of the implication
% $(i)\implies(ii)$ in \cite[Theorem~1.2]{HR}.
% More generally, this result applies to every unital
% $C^*$-algebra in which each element admits a polar
% decomposition within the algebra.

% It is also known that, for a compact Hausdorff space
% $\Omega$, every $f\in C(\Omega)$ admits a factorization
% $f=u|f|$ for some $u\in C(\Omega)$ if and only if
% $\Omega$ is an $F$-space; see
% \cite[Theorem~2.14]{ASCG}. However, the function $u$
% in this characterization is not required to be a
% partial isometry. Thus, the characterization does not
% directly imply the hypothesis of Theorem~\ref{thm:A-has-p.d.}.
% Nevertheless, $C^*$-extremality and linear extremality
% coincide in the closed unit ball of $C(\Omega)$ for
% every compact Hausdorff space $\Omega$, without
% requiring the $F$-space property.
% \end{remark}
% \end{remark}

The following corollary extends the preceding result to quotient \(C^*\)-algebras. It shows that if every element \(x\in\cA\) admits a factorization \(x=v|x|\) for some partial isometry \(v\in\cA\), then \(C^*\)-extremality implies linear extremality in every unital quotient of \(\cA\).
\begin{corollary}
\label{cor:ideal}
    Let $\cA$ be a $C^*$-algebra such that every element $x \in \cA$ has a factorization $x=v |x|$, where $v\in \cA$ is a partial isometry. Let $J\subset \cA$ be a closed two-sided ideal in $\cA$. Then every $C^*$-extreme point of $(\cA/J)_1$ is a linear extreme point of $(\cA/J)_1$. 
\end{corollary}

\begin{proof}
% \begin{proof}
Let $\pi:\A\to\A/J$ be the quotient map, and let
$a\in(\A/J)_1$ be a $C^*$-extreme point.
Choose $x\in\A$ such that $\pi(x)=a$.
By hypothesis, there exists a partial isometry
$v\in\A$ such that $x=v|x|$.

Since $\pi$ is a $*$-homomorphism,
$\pi(v)$ is a partial isometry and
\[
\pi(|x|)=|\pi(x)|=|a|.
\]
Consequently,
\[
a=\pi(x)=\pi(v)\pi(|x|)=\pi(v)|a|.
\]
Thus, $a$ admits the required factorization
in $\A/J$. Applying
Theorem~\ref{thm:A-has-p.d.}, we conclude
that $a$ is a linear extreme point of
$(\A/J)_1$.
% \end{proof}
    % Let $a \in (\cA/J)_1$ be a $C^*$-extreme point of $(\cA/J)_1$. Let $x$ be a lift of $a$ in $\cA$. Then by hypothesis, $x=v|x|$ for some partial isometry $v \in \cA$. Let $\pi: \cA \to \cA/J$ be the quotient map. Then $\pi(v)$ is a partial isometry in $\cA/J$ and $\pi(|x|)=|a|$. Now, using the fact that every positive contraction in a unital $C^*$-algebra can be written as an average of two unitaries, one has $|a|=\frac{1}{2}(u_1+u_2)$, where $u_1, u_2 \in \cA/J$ are unitaries. Therefore, 
    % \[
    % a=\frac{1}{2}(\pi(v)u_1+\pi(v)u_2),
    % \]
    % which is a proper $C^*$-convex decomposition of $a$. Applying the argument as in Theorem \ref{thm:A-has-p.d.}, one obtains that $a$ is a partial isometry in $\cA/J$. Now using Proposition \ref{p:p.i.+c_ext}, we get $a$ is a linear extreme point of $\cA/J$. This completes the proof.
\end{proof}

\begin{corollary}
    \label{cor:cal-alg}
    Let $\mathcal{Q}(H)$ denote the Calkin algebra. Let $x\in \cQ(H)_1$ be a $C^*$-extreme point of $\cQ(H)_1$. Then $x$ is a linear extreme point of $\cQ(H)_1$.
\end{corollary}

\begin{proof}
    It follows directly from Corollary \ref{cor:ideal} and the observation that in $B(H)$ every element has the polar decomposition.
\end{proof}

We next consider \(C^*\)-subalgebras of \(B(H)\) containing the compact operators. In this setting, Theorem \ref{thm:lin-ext} takes a particularly simple form. 
% Indeed, the presence of rank-one operators forces at least one of the defect projections of a linear extreme point to vanish. This leads to the following characterization.
\begin{theorem}
\label{thm:lin-ext-char}
Let $\cA\subseteq B(H)$ be a unital $C^*$-algebra
such that $\cK(H)\subseteq\cA$.
Then the linear extreme points of $\cA_1$
are precisely the isometries and coisometries
in $\cA$.
\end{theorem}

\begin{proof}
By Theorem~\ref{thm:lin-ext}, every isometry
and every coisometry in $\cA$ is a linear
extreme point of $\cA_1$.
Conversely, let $x$ be a linear extreme
point of $\cA_1$. By Theorem~\ref{thm:lin-ext},
$x$ is a partial isometry satisfying Equation \eqref{eq:ext-cond}.

Suppose, to the contrary, that both
$\mathbf{1}_{\cA}-xx^*$ and $\mathbf{1}_{\cA}-x^*x$
are nonzero. Since $x$ is a partial isometry,
these are nonzero projections in $\cA$.
Choose nonzero vectors
\[
\xi\in(\mathbf{1}_{\cA}-xx^*)H,
\qquad
\eta\in(\mathbf{1}_{\cA}-x^*x)H.
\]
Consider the rank-one operator
$\xi\eta^*\in\cK(H)\subseteq\cA$.
Then
\[
\begin{aligned}
&(\mathbf{1}_{\cA}-xx^*)\xi\eta^*
(\mathbf{1}_{\cA}-x^*x)\eta \\
&\qquad
=(\mathbf{1}_{\cA}-xx^*)\xi\eta^*\eta \\
&\qquad
=\|\eta\|^2\xi\neq0.
\end{aligned}
\]
This contradicts the above condition.
Therefore, either
\[
\mathbf{1}_{\cA}-xx^*=0
\quad\text{or}\quad
\mathbf{1}_{\cA}-x^*x=0.
\]
Hence $x$ is either a coisometry or an
isometry.
\end{proof}

\begin{corollary}
Let $\cT$ be the Toeplitz $C^*$-algebra.
Then the linear extreme points of $\cT_1$
are precisely the isometries and coisometries
in $\cT$.
\end{corollary}

\begin{proof}
By \cite[Proposition~11.2.19]{E},
the Toeplitz $C^*$-algebra $\cT$ contains
$\cK(H^2)$.
Hence the conclusion follows directly
from Theorem~\ref{thm:lin-ext-char}.
\end{proof}

\vspace{.3cm}
 % \noindent
\textbf{Acknowledgements:} This work is part of the project ``Classification of Banach spaces using differentiability",
funded by the Anusandhan National Research Foundation (ANRF), Core Research Grant, CRG2023-000595.
The author is a research associate (RA) in this project. The author is grateful to Prof. T. S. S. R. K. Rao (SNU) for his valuable insights and numerous helpful discussions throughout the preparation of this article.

\vspace{.3cm}
\textbf{Declaration:}
No conflicts of interest.

\bibliographystyle{alpha}
\bibliography{references}
\end{document}